\documentclass[12pt]{article}
\usepackage{color}
\definecolor{darkblue}{rgb}{0.00,0.25,0.50}
\usepackage[colorlinks,filecolor=blue,citecolor=darkblue]{hyperref}

\usepackage{amsfonts,amssymb,amsmath,amsthm}
\usepackage{url}
\usepackage{enumerate}
\usepackage[ukrainian, russian, english]{babel}
\usepackage[cp1251]{inputenc}
\newtheorem{theorem}{Теорема}
\newtheorem{lem}{Лемма}
\newtheorem{cor}{Следствие}

\begin{document}\selectlanguage{russian}
\thispagestyle{empty}

\title{}

\begin{center}
\textbf{\Large Равномерное приближение периодических функций тригонометрическими суммами
специального вида}\footnote{Работа частично поддержана Государственным фондом
фундаментальных исследований Украины (проект № GP/$\Phi$36/068)}
\end{center}
\vskip0.5cm
\begin{center}
А.С. Сердюк, Е.Ю. Овсий\\ \emph{\small Институт математики НАН Украины, Киев}
\end{center}
\vskip0.5cm


\begin{abstract}
В работе изучаются аппроксимационные свойства тригонометрических сумм $U_{n,p}^\psi $
специального вида на классах $(\psi ,\beta )$-дифференцируемых (в смысле Степанца)
периодических функций $C^{\psi}_{\beta,\infty}.$ Вследствие согласованности между
параметрами аппроксимационных сумм и приближаемых классов в достаточно общей ситуации
удается найти решение соответствующей задачи Колмогорова--Никольского. Показано, что в ряде
важных случаев рассматриваемые суммы на классах $C^{\psi}_{\beta,\infty}$ обеспечивают более
высокий порядок приближения в метрике пространства $C$ по сравнению с суммами Фурье,
Зигмунда и Валле Пуссена. Указан диапазон параметров, в пределах которого суммы
$U_{n,p}^\psi $ доставляют порядок наилучшего равномерного приближения классов
$C^{\psi}_{\beta,\infty}$.
\end{abstract}


\section{Введение и постановка задачи}

Пусть $C$ --- пространство непрерывных $2\pi $-периодических функций $f$, в котором норма
определяется формулой $\|f\|_C=\max\limits_{x}{|f(x)|}.$

Рассмотрим класс $C^{\psi}_{\beta,\infty}$ \cite{Stepanets_1987_Izv} непрерывных $2\pi
$-периодических функций $f\in C$, для которых при $\beta \in\mathbb{R}$ и заданной
последовательности $\psi (k)$ ($k\in\mathbb{N}$) действительных чисел, ряд
    $$\sum\limits_{k=1}^{\infty}\frac{1}{\psi (k)}\bigg(
    a_k(f)\cos \Big(kx+\frac{\beta \pi }{2}\Big)+
    b_k\sin \Big(kx+\frac{\beta \pi }{2}\Big)\bigg)$$
является рядом Фурье некоторой функции $\varphi \in S_\infty$, где
    $$S_\infty=\{\varphi :\ \mathop{\rm ess\,sup}\limits_{t\ }|\varphi (t)|\leqslant 1\}.$$
Функцию $\varphi $ принято называть $(\psi ,\beta)$-производной функции $f$ и обозначать
через $f^\psi _\beta .$ При $\psi (k)=k^{-r},$ $r>0,$ класс $C^{\psi}_{\beta,\infty}$
совпадает с классом Вейля-Надя $W^r_\beta $, а при $\beta =r$ --- с классом Вейля $W^r_r.$ В
случае натуральных $r$ и $\beta =r$, класс $C^{\psi}_{\beta,\infty}$ является классом $W^r$
периодических функций, чьи $r$-е производные по модулю почти всюду не превышают единицы.
Если
    \begin{equation}\label{27.11.12-23:22:37}
    \lim\limits_{k\to\infty}{k^r\psi (k)}=0\ \ \ \forall r\in \mathbb{R},
    \end{equation}
то класс $C^{\psi}_{\beta,\infty}$ состоит (см. \cite[Гл. 1, Разд. 8]{Stepanets_1995}) из
бесконечно дифференцируемых функций. Примером последовательности $\psi (k)$, удовлетворяющей
условию (\ref{27.11.12-23:22:37}) является последовательность $\psi (k)=e^{-\alpha k^r },$
$\alpha
>0$, $r>0$. В этом случае классы $C^{\psi}_{\beta,\infty}$ будем обозначать через $C^{\alpha
,r}_{\beta ,\infty}.$ Если $\psi (k)$ удовлетворяет условию
    \begin{equation}\label{27.11.12-23:22:27}
    |\psi (k)|\leqslant Ke^{-\alpha k},\ \ \ k\in\mathbb{N},\ \ \ \alpha >0,
    \end{equation}
то $C^{\psi}_{\beta,\infty}$ состоит из аналитических функций, регулярно продолжающихся в
полосу $|\text{Im} z|<\alpha $ комплексной плоскости.

Следуя работе \cite[c. 147]{Stepanets_VSP}, обозначим через $\mathfrak{M}$ множество всех
положительных, выпуклых вниз функций $\psi (t),$ $t\geqslant 1$, удовлетворяющих условию
    $$\lim\limits_{t\to\infty}{\psi (t)}=0$$
и поставим в соответствие каждой $\psi \in\mathfrak{M}$ характеристику вида
    \begin{equation}\label{27.11.12-23:55:28}
    \eta(t)=\eta(\psi ;t)=\psi ^{-1}\bigg(\frac{1}{2}\psi (t)\bigg),\ \ \ t\geqslant 1,
    \end{equation}
где $\psi ^{-1}(\cdot)$ --- функция, обратная к $\psi (\cdot).$ При помощи характеристики
$\eta(t)$ выделим из $\mathfrak{M}$ подмножество $F$ следующим образом:
    $$F=\{\psi \in \mathfrak{M}:\ \ \eta'(\psi ;t)\leqslant K,\ \ \forall t\geqslant 1\}.$$
Как показано в работе \cite[c. 153]{Stepanets_VSP}, в $F$ входят все функции $\psi
\in\mathfrak{M}$, для которых
    \begin{equation}\label{27.11.12-23:52:37}
    0<C_1\leqslant \frac{t}{T(t)}\leqslant C_2,\ \ \forall t\geqslant 1,\ \ C_1,C_2
    \text{ = const},
    \end{equation}
где $T(t)=T(\psi ;t)=\eta(\psi ;t)-t$. Множество таких функций обозначают через
$\mathfrak{M}_C.$ Величина $T(t)$ имеет простую геометрическую интерпретацию, она равна
длине промежутка, на котором значение функции $\psi (t)$ уменьшается ровно в два раза, в
связи с этим, величину $T(t)$ естественно назвать периодом полураспада функции $\psi$.
Примерами $\psi (\cdot)$, принадлежащих $\mathfrak{M}_C$, могуть служить, например, функции
вида $\psi _1(t) = t^{-r},$ $r>0$, $\psi _2(t)=\frac{1}{t^{r}\ln(t+\beta )}$, $r>0,$ $\beta
\geqslant 1$ и другие. Множество $F$ включает также (см. \cite[с. 153]{Stepanets_VSP})
подмножество $\mathfrak{M}_\infty^+$ всех функций $\psi \in \mathfrak{M}$, для которых
характеристика
    $$\mu (t)=\mu (\psi ;t) =\frac{t}{\eta(t)-t},\ \ t\geqslant 1,$$
именуемая модулем полураспада, монотонно стремиться к бесконечности при $t\to\infty.$ Для
$\psi _3(t)=e^{-\alpha t^r},$ $\alpha >0,$ $r>0,$ $\eta(\psi _3;t)=t\Big(\frac{\ln 2}{\alpha
t^r}+1\Big)^{1/r}$ и тогда при $t\to\infty$
    $$\mu (\psi _3;t)=\frac{t}{\eta(\psi _3;t)-t}
    =\frac{1}{(\frac{\ln 2}{\alpha t^r}+1)^{1/r}-1}\uparrow\infty.$$
Таким образом $\psi _3\in\mathfrak{M}_\infty^+\subset F.$ Следовательно, среди функций $\psi
$, принадлежащих множеству $F,$ находятся функции, имеющие степенную скорость стремления к
нулю, а также те, которые стремятся к нулю быстрее любой степенной функции. Пример функции
$\psi _4(t)=\frac{1}{\ln(t+1)}$, для которой $\eta(\psi _4;t)=(t+1)^2-1$, а значит
$\eta'(\psi _4;t)=2(t+1)$, показывает, что множество $F$ может не содержать функций,
стремящихся к нулю медленнее любой степенной функции.

В дальнейшем, не уменьшая общности, будем считать, что последовательность $\psi (k)$,
задающая класс $C^{\psi}_{\beta,\infty}$, является сужением на множестве натуральных чисел
$\mathbb{N}$ некоторой функции $\psi (t),$ $t\geqslant 1,$ из множества $F.$

Рассмотрим для произвольной функции $f\in C$ сумму вида
    \begin{equation}\label{8.09.12-14:41:22}
        U_{n,p}^{\varphi }(f;x)=\sum\limits_{
        k=0}^{n-1}\lambda _{n,p}(k)A_k(f;x),
    \end{equation}
где
    $$\lambda _{n,p}(k)=\lambda _{n,p}(k;\varphi  )=
  \begin{cases}
    1, & 0\leqslant k\leqslant n-p, \\
    1-\frac{\varphi  (k)}{\varphi  (n)}, & n-p+1\leqslant k\leqslant n-1,
  \end{cases}
    $$
$\varphi  (k)$ $(k\in\mathbb{N})$ --- произвольная монотонно возрастающая к бесконечности
последовательность действительных чисел $p\in\mathbb{N},$ $p\in[1,n],$
    $$A_k(f;x):=a_k(f)\cos kx+b_k(f)\sin kx,$$
    $$A_0(f;x):=\frac{a_0(f)}{2}$$
и $a_0(f),$ $a_k(f)$ и $b_k(f)$ --- коэффициенты Фурье функции $f.$

Суммы $U_{n,p}^\varphi (f;x)$ при определенном выборе параметров $p$ и $\varphi  (k)$
совпадают с такими классическими суммами, как суммы Зигмунда \cite{Zygmund_1945} (при $p=n$
и $\varphi  (k)=k^{s},$ $s>0$)
    $$Z_n^s(f;x)=\sum\limits_{k=0}^{n-1}\bigg(1-\frac{k^s}{n^s}\bigg)A_k(f;x),\ \ \ s>0,$$
суммы Фейера \cite{Fejer_1903} (при $p=n$ и $\varphi  (k)=k$)
    $$\sigma_{n-1}(f;x)=\sum\limits_{k=0}^{n-1}\bigg(1-\frac{k}{n}\bigg)A_k(f;x),$$
суммы Валле Пуссена \cite{Vallee Poussin} (при $p\in\mathbb{N},$ $1\leqslant p\leqslant n$ и
$\varphi (k)=k-n+p$)
    $$V_{n,p}(f;x)=\sum\limits_{k=0}^{n-1}\lambda _{n,p}(k)A_k(f;x),$$
где
    $$\lambda _{n,p}(k)=
  \begin{cases}
    1, & 0\leqslant k\leqslant n-p, \\
    1-\frac{k-n+p}{p}, & n-p+1\leqslant k\leqslant n-1,
  \end{cases}
    $$
суммы Фурье (при $p=1$)
    $$S_{n-1}(f;x)=\sum\limits_{k=0}^{n-1}A_k(f;x).$$
При $p=n$ суммы $U_{n,p}^\varphi (f;x)$ совпадают с так называемыми обобщенными суммами
Зигмунда \cite{Aljancic_1959} (см., также, \cite{Gavrilyuk, Serdyuk_Ovsii_2009})
    $$Z_n^\varphi  (f;x)=\sum\limits_{k=0}^{n-1}\bigg(
    1-\frac{\varphi  (k)}{\varphi  (n)}\bigg)A_k(f;x),$$
где $\varphi  (k)$ ($k\in\mathbb{N}$) --- произвольная монотонно возрастающая к
бесконечности последовательность действительных чисел.

Целью данной работы является изучение асимптотического поведения при $n\to\infty$ величины
\begin{equation}\label{23.11.12-00:12:21}
    \mathcal{E}(C^{\psi}_{\beta,\infty};U_{n,p}^\psi)=\sup\limits_{
    f\in C^{\psi}_{\beta,\infty}}{\|f(\cdot)-U_{n,p}^\psi(f;\cdot)\|_C},\ \ \psi \in F, \ \ \beta
    \in\mathbb{R}
\end{equation}
где $U_{n,p}^{\psi }(f;\cdot)$ --- суммы $U_{n,p}^{\varphi }(f;\cdot)$ вида
(\ref{8.09.12-14:41:22}) при $\varphi (k)=\frac{k-n+p}{\psi (k)}.$ Отметим, что данная
тематика для сумм Фурье, Валле Пуссена и Зигмунда на различных функциональных классах имеет
большую историю, связанную с именами \mbox{А.Н. Колмогорова,} \mbox{С.М. Никольского,}
\mbox{А.Ф. Тимана,} \mbox{В.К. Дзядыка,} \mbox{С.Б. Стечкина,} \mbox{Н.П. Корнейчука,}
\mbox{А.В. Ефимова,}  \mbox{С.А. Теляковского,} \mbox{А.И. Степанца,} \mbox{В.П. Моторного,}
\mbox{Р.М. Тригуба,} \mbox{В.И. Рукасова} и многих других. Детально с историей данного
вопроса можно ознакомиться, в частности, по работам \cite{Stepanets_1999, Stepanets_2001,
Stepanets_2002} и \cite{STEANETS-RUKASOV}.

Если для величины (\ref{23.11.12-00:12:21}) получено асимптотическое равенство, то есть
равенство вида
    $$\mathcal{E}(C^{\psi}_{\beta,\infty};U_{n,p}^\psi)=
    \nu (n)+o(1)\nu (n),\ \ \ n\to\infty,$$
где $\nu (n)=\nu (n,p,\psi ,\beta )$ некая конкретная последовательность, то следуя А.И.
Степанцу \cite{Stepanets_1987_Izv} будем говорить, что для сумм $ U_{n,p}^\psi(f;x)$ найдено
решение задачи Колмогорова--Никольского на классе $C^{\psi}_{\beta,\infty}.$

Перейдем к изложению основных результатов.

\section{Основные результаты}

Имеют место следующие утверждения.

\begin{theorem}\label{17.11.12-11:21:14} Пусть
$\psi \in F,$ $\beta\in \mathbb{R},$ $n,p\in \mathbb{ N},$ $p\leqslant n.$ Тогда при
$n\to\infty$
\begin{equation}\label{14.11.12-22:35:31}
    \mathcal{E}(C^\psi _{\beta ,\infty};U_{n,p}^\psi)=\psi (n)\bigg(\frac{4}{\pi ^2}
    A^\psi _{n,p}+O(1)\bigg),
\end{equation}
где
\begin{equation}\label{17.11.12-11:22:23}
    A^\psi _{n,p}=
  \begin{cases}
    \ln p, &  \text{если}\,\,\,\, T(n)\leqslant   1,\\
    &\\
    \ln\frac{p}{T(n)}, & \text{если}\,\,\,\, 1\leqslant T(n)\leqslant  p,\\
    &\\
    \ln\frac{T(n)}{p}, & \text{если}\,\,\,\, T(n)\geqslant p,
  \end{cases}
\end{equation}
$T(t)=\eta(\psi ;t)-t,$ $\eta(\psi ;n)=\psi ^{-1}\big(\frac{1}{2}\psi (n)\big),$ $O(1)$ ---
величина, равномерно ограниченная по $\beta ,$ $n$ и $p.$
\end{theorem}

Как показано в работе \cite[c. 508]{Stepanets_VSP}, если $\psi \in F$ и $\beta
\in\mathbb{R},$ то для величины
    $$E_n(C^{\psi}_{\beta,\infty})=\sup\limits_{f\in C^{\psi}_{\beta,\infty}}{\inf_{
    t_{n-1}\in \mathcal{T}_{n-1}}\|f(\cdot)-t_{n-1}(\cdot)\|_C}$$
наилучшего равномерного приближения класса $C^{\psi}_{\beta,\infty}$ тригонометрическими
полиномами, порядок которых не превышает $n-1$, имеет место порядковая оценка
\begin{equation}\label{12.12.12-00:55:36}
    E_n(C^{\psi}_{\beta,\infty})\asymp\psi (n)
\end{equation}
(запись $\alpha (n)\asymp\beta (n)$ означает, что существуют константы $K_1,K_2>0$ такие,
что $K_1\beta (n)\leqslant \alpha (n)\leqslant K_2\beta (n)$). Исходя из теоремы
\ref{17.11.12-11:21:14} и оценки (\ref{12.12.12-00:55:36}), приходим к следующему
утверждению.

\begin{cor}\label{12.12.12-00:23:46} Пусть
$\beta\in \mathbb{R},$ $n,p\in \mathbb{ N},$ $p\leqslant n$. Тогда если $\psi \in F$ и
$p=p(n)$ такова, что $T(n)\asymp p(n)$, то
    $$\mathcal{E}(C^\psi _{\beta ,\infty};U_{n,p}^\psi)\asymp E_{n}(C^{\psi}_{\beta,\infty})\asymp \psi (n),$$
то есть суммы $U_{n,p}^\psi (f;x)$ реализуют порядок наилучшего равномерного приближения
класса $C^\psi _{\beta ,\infty}$.
\end{cor}

При $n\to\infty$, $p\to\infty$ и $T(n)=o(1)p(n)$ равенство (\ref{14.11.12-22:35:31}) дает
решение задачи Колмогорова--Никольского для сумм $ U_{n,p}^\psi(f;x)$, поскольку в этом
случае
    $$\lim\limits_{n\to\infty}{A_{n,p}^\psi }=\lim\limits_{
    p\to\infty}{\ln p}=\infty.$$


Заметим, что на классе $C^{\alpha ,1}_{\beta ,\infty}$ при $p\to\infty$ и $n-p\to\infty$
суммы $U_{n,p}^{\psi}(f;x)=U_{n,p}^{\alpha ,1}(f;x)$ обеспечивают лучший порядок приближения
по сравнению с классическими суммами Валле Пуссена $V_{n,p}(f;x)$. Действительно, для сумм
$V_{n,p}(f;x)$ имеет место, в частности, следующее асимптотическое равенство (см., например,
\cite[c. 130]{Serdyuk V-P}, \cite[c. 10]{Rendiconti})
\begin{equation}\label{17.11.12-16:02:15}
    \mathcal{E}(C^{\alpha ,1 } _{\beta ,\infty};V_{n,p})=
    \frac{e^{-\alpha (n-p+1)}}{p}\bigg(\frac{4}{\pi (1-e^{-2\alpha })}+$$
    $$+O(1)
    \Big(\frac{e^{-\alpha }}{(1-e^{-\alpha })^3(n-p+1)}+
    \frac{e^{-\alpha p}}{1-e^{-\alpha }}\Big)\bigg),
\end{equation}
где $O(1)$ --- величина, равномерно ограниченная по $n,$ $p,$ $\alpha $ и $\beta .$
Сопоставляя (\ref{14.11.12-22:35:31}) и (\ref{17.11.12-16:02:15}) находим, что если $p=p(n)$
удовлетворяет условиям
    $$p\to\infty,\ \ \ n-p\to\infty,$$
то
    $$\lim\limits_{\stackrel{{n-p\to\infty}}{p\to\infty}}\frac{\mathcal{E}(C^{\alpha,1
     } _{\beta ,\infty};U_{n,p}^{\alpha,1 })}
    {\mathcal{E}(C^{\alpha,1 } _{\beta ,\infty};V_{n,p})}=0.$$
Как было отмечено ранее, при $p=1$ суммы $U_{n,p}^\psi(f;x)$ совпадают с суммами Фурье
$S_{n-1}(f;x)$ порядка $n-1.$ При указанном значении параметра $p$, из теоремы
\ref{17.11.12-11:21:14} получаем следующее утверждение.

\begin{cor}\label{17.11.12-17:02:06} Пусть
$\psi \in F,$ $\beta\in \mathbb{R}$ и $n\in \mathbb{ N}.$ Тогда при $n\to\infty$
\begin{equation}\label{17.11.12-17:04:30}
    \mathcal{E}(C^\psi _{\beta ,\infty};S_{n-1})=\psi (n)\bigg(\frac{4}{\pi ^2}
    \ln^{+}(\eta(\psi ;n)-n)+O(1)\bigg),
\end{equation}
где $\ln^+(t)=\max\{\ln t, 0\},$ $\eta(\psi ;n)=\psi ^{-1}\big(\frac{1}{2}\psi (n)\big),$ а
$O(1)$
--- величина, равномерно ограниченная по $\beta $ и $n.$
\end{cor}

Равенство (\ref{17.11.12-17:04:30}) получено А.И. Степанцом (см., например, \cite[c.
257]{Stepanets_VSP}) и дает решение задачи Колмогорова--Никольского для сумм Фурье в случае,
когда при $n\to\infty$
    $$\eta(\psi ;n)-n\to\infty.$$
При $\psi (t)=t^{-r},$ $\beta =r$ (в этом случае $C^{\psi}_{\beta,\infty}=W^r_r$, $\eta(\psi
;n)-n=(2^{1/r}-1)n$) равенство (\ref{17.11.12-17:04:30}),  принимает вид
    $$\mathcal{E}(W^r; S_{n-1})=n^{-r}\bigg(
    \frac{4}{\pi ^2}\ln n+O(1)\bigg),$$
и установлено А.Н. Колмогоровым \cite{Kolmogoroff_1935} (при $r\in\mathbb{N}$) и В.Т.
Пинкевичем \cite{Pinkewitch} (при $r>0$).

Сопоставление следствия \ref{12.12.12-00:23:46} с равенством (\ref{17.11.12-17:04:30})
показывает, что в случае, когда $T(n)\to\infty,$ а параметр $p=p(n)$ выбран таким образом,
чтобы $p(n)\asymp T(n)$, суммы $U_{n,p}^\psi (f;x)$ осуществляют лучший порядок приближения
на классе $C^{\psi}_{\beta,\infty},$ чем суммы Фурье.

Положив в теореме \ref{17.11.12-11:21:14} $p=n$ и учитывая, что в этом случае $U^\psi
_{n,n}(f;x)$
--- обобщенные суммы Зигмунда $Z_n^\varphi   (f;x)$ $(\varphi   (k)=k/\psi (k),$
$k\in\mathbb{N}),$ получаем.

\begin{cor}\label{17.11.12-23:50:15} Пусть
$\psi \in F,$ $\varphi  (k)=k/\psi (k),$ $\beta\in \mathbb{R}$ и $n\in \mathbb{ N}.$ Тогда
при $n\to\infty$
\begin{equation}\label{17.11.12-23:50:18}
    \mathcal{E}(C^\psi _{\beta ,\infty};Z^\varphi   _{n})=\psi (n)\bigg(\frac{4}{\pi ^2}
    A_n^\psi +O(1)\bigg),
\end{equation}
где
$$
    A_{n}^\psi =
  \begin{cases}
    \ln n, &  \text{если}\,\,\,\, T(n)\leqslant   1,\\
    &\\
    \ln\frac{n}{T(n)}, & \text{если}\,\,\,\, 1\leqslant T(n)\leqslant  n,\\
    &\\
    \ln\frac{T(n)}{n}, & \text{если}\,\,\,\, T(n)\geqslant n,
  \end{cases}
$$
$T(t)=\eta(\psi ;t)-t,$ $\eta(\psi ;n)=\psi ^{-1}\big(\frac{1}{2}\psi (n)\big),$ а $O(1)$
--- величина, равномерно ограниченная по $\beta $ и $n.$
\end{cor}

Поскольку в силу определения (\ref{27.11.12-23:52:37}), для любой функции $\psi \in
\mathfrak{M}_C$
    $$T(n)\asymp n$$
и $\mathfrak{M}_C\subset F,$ то из (\ref{17.11.12-23:50:18}) находим
\begin{equation}\label{12.12.12-23:19:08}
    \mathcal{E}(C^\psi _{\beta ,\infty};Z^\varphi   _{n})=O(1)\psi (n),\ \ \ \psi \in
    \mathfrak{M}_C.
\end{equation}
Впрочем, эта же оценка вытекает также из следствия \ref{12.12.12-00:23:46} при $p=n.$

Положив в (\ref{12.12.12-23:19:08}) $\psi (t)=t^{-r},$ $r>0$ (в этом случае $
C^{\psi}_{\beta,\infty}=W^{r}_{\beta ,\infty}$ и, соответственно, $Z_n^\varphi   (f;x)=
Z_n^s(f;x)$, $s=r+1$), получаем
\begin{equation}\label{18.11.12-00:33:06}
    \mathcal{E}(W^{r}_{\beta ,\infty};Z^s _{n})=O(1)n^{-r}.
\end{equation}
Отметим, что оценка (\ref{18.11.12-00:33:06}) может быть получена из более точных
результатов \mbox{С.А. Теляковского} \cite{TELYKOVSKIY-61} (см. теорему 7).

\section{Доказательство теоремы \ref{17.11.12-11:21:14}}

Пусть $f\in C^\psi _{\beta,\infty} $ и $\psi \in \mathfrak{M}',$ где
    $$\mathfrak{M}'=\bigg\{\psi \in\mathfrak{M}:\ \ \int\limits_{1}^{\infty}\frac{\psi (t)}{t}\,dt
    <\infty\bigg\}$$
(как показано в \cite[c. 155]{Stepanets_VSP}, $F\subset\mathfrak{M}'$). Тогда, в силу
теоремы 4.1 работы \cite[c. 71]{Stepanets_1995}, в каждой точке $x\in \mathbb{R}$ для
величины
    $$\rho _{n,p}(f;x):=f(x)-U_{n,p}^\psi(f;x)$$
имеет место равенство
\begin{equation}\label{8.09.12-15:04:47}
    \rho _{n,p}(f;x)=\int\limits_{-\infty}^{\infty}
    f^\psi _\beta \Big(x+\frac{t}{n}\Big)\widehat{\tau}_{n,p}(t)\,dt, \ \ n\in\mathbb{N},
\end{equation}
где
    $$\widehat{\tau}_{n,p}(t):=\frac{1}{\pi }\int\limits_{0}^{\infty}
    \tau _{n,p}(u)\cos\Big(ut+\frac{\beta \pi }{2}\Big)\,du,$$
а
\begin{equation}\label{2.12.12-20:32:42}
    \tau _{n,p}(u)=
  \begin{cases}
    0, & 0\leqslant u\leqslant \frac{n-p}{n}, \\
    \psi (n)\frac{nu-n+p}{p}, & \frac{n-p}{n}\leqslant u\leqslant 1,\\
    \psi (nu), & u\geqslant 1.
  \end{cases}
\end{equation}
Упростим правую часть равенства (\ref{8.09.12-15:04:47}) с целью выделения главного члена
величины $\rho _{n,p}(f;x).$ С этой целью, положим
\begin{equation}\label{8.09.12-15:36:59}
    \widehat{\tau}_{n,p+}(t):=\frac{1}{\pi }\int\limits_{
    0}^{\infty}\tau _{n,p}(u)\cos ut\,du
\end{equation}
и
\begin{equation}\label{8.09.12-15:38:11}
    \widehat{\tau}_{n,p-}(t):=\frac{1}{\pi }\int\limits_{
    0}^{\infty}\tau _{n,p}(u)\sin ut\,du.
\end{equation}
С учетом (\ref{8.09.12-15:36:59}) и (\ref{8.09.12-15:38:11}), равенство
(\ref{8.09.12-15:04:47}) можно представить следующим образом
\begin{equation}\label{8.09.12-15:40:50}
    \rho _{n,p}(f;x)=\cos\frac{\beta \pi }{2}\int\limits_{
    -\infty}^{\infty}f^\psi _\beta \Big(x+\frac{t}{n}\Big)\widehat{\tau}_{n,p+}(t)\,dt-$$
    $$-\sin\frac{\beta \pi }{2}\int\limits_{
    -\infty}^{\infty}f^\psi _\beta \Big(x+\frac{t}{n}\Big)\widehat{\tau}_{n,p-}(t)\,dt=$$
    $$=
    \cos\frac{\beta \pi }{2}\rho _{n,p+}(f;x)-\sin\frac{\beta \pi }{2}\rho _{n,p-}(f;x).
\end{equation}

Поскольку класс $C^{\psi}_{\beta,\infty}$ инвариантен относительно сдвига аргумента (если
$f\in C^{\psi}_{\beta,\infty}$, то и функция $f_1(\cdot)= f(\cdot+h),$ $h\in\mathbb{R},$
также принадлежит классу $C^{\psi}_{\beta,\infty}$), то для величины
(\ref{23.11.12-00:12:21}) можно записать
\begin{equation}\label{1.12.12-16:19:29}
    \mathcal{E}(C^{\psi}_{\beta,\infty};U_{n,p}^\psi )=
    \sup\limits_{f\in C^{\psi}_{\beta,\infty}}{|\rho _{n,p}(f;0)|}.
\end{equation}
Исходя из (\ref{1.12.12-16:19:29}), достаточно ограничиться рассмотрением отклонения $\rho
_{n,p}(f;x)$ в точке $x=0.$

Докажем следующее утверждение.

\begin{lem}\label{8.09.12-15:54:03} Пусть $\psi \in \mathfrak{M}',$ $\beta\in \mathbb{R},$
$n,p\in \mathbb{ N},$ $p\leqslant n$ и $\alpha  (n)$ --- произвольная последовательность
действительных чисел, удовлетворяющих условию
    $$\alpha (n)\geqslant \frac{K}{n},\ \ n\in \mathbb{N},$$
где $K$ --- некоторая положительная константа. Тогда для произвольной функции $f\in C^\psi
_{\beta,\infty}$ при $n\to\infty$ имеет место равенство:
\begin{equation}\label{26.09.12-13:08:00}
   \rho _{n,p}(f;0)=(-1)^{s}\frac{\psi (n)}{\pi }\int\limits_{\mathcal{I}_{\alpha_n}}
    f^\psi _\beta \Big(\frac{t}{n}\Big)\frac{\sin \Big(t+\frac{\beta\pi}{2}\Big)}{t}\,dt+O(1)r_n,
\end{equation}
где
    $$\mathcal{I}_{\alpha_n}=\Big(-\frac{n}{p},-\alpha (n)n\Big)\bigcup
    \Big(\alpha (n)n, \frac{n}{p}\Big),\ \ \text{если $\alpha (n)\leqslant \frac{1}{p}$},$$
    $$\mathcal{I}_{\alpha_n}=\Big(-\alpha (n)n, -\frac{n}{p}\Big)\bigcup
    \Big(\frac{n}{p}, \alpha (n)n\Big),\ \ \text{если $\frac{1}{p}<\alpha (n)\leqslant 1 $},$$
    $$\mathcal{I}_{\alpha_n}=\Big(-n, -\frac{n}{p}\Big)\bigcup
    \Big(\frac{n}{p}, n\Big),\ \ \text{если $\alpha (n)>1 $},$$
    \begin{equation}\label{8.09.12-15:12:50}
    s=s(\alpha ,p)=
  \begin{cases}
    1, & \alpha (n)\leqslant \frac{1}{p}, \\
    0, & \alpha (n)> \frac{1}{p},
  \end{cases}
\end{equation}
\begin{equation}\label{8.09.12-16:24:42}
    r_n=\psi (n)+\int\limits_{
    1/\alpha (n)}^{\infty}\frac{\psi (t+n)}{t}\,dt+\int\limits_{
    \alpha (n)}^{\infty}\frac{\psi (n)-\psi (n+1/t)}{t}\,dt,
\end{equation}
а $O(1)$ --- величина, равномерно ограниченная относительно $\beta$, $n$ и $p.$
\end{lem}

Перейдем к доказательству леммы.

\begin{proof} В силу (\ref{8.09.12-15:40:50}), лемма будет доказана, если установить оценки
\begin{equation}\label{8.09.12-16:00:50}
    \rho _{n,p+}(f;0)=(-1)^s\frac{\psi (n)}{\pi }\int\limits_{\mathcal{I}_{\alpha_n}}
    f^\psi _\beta \Big( \frac{t}{n}\Big)\frac{\sin t}{t}\,dt+O(1)r_n,
\end{equation}
\begin{equation}\label{8.09.12-16:01:00}
    \rho _{n,p-}(f;0)=(-1)^{s+1}\frac{\psi (n)}{\pi }\int\limits_{\mathcal{I}_{\alpha_n}}
    f^\psi _\beta \Big( \frac{t}{n}\Big)\frac{\cos t}{t}\,dt+O(1)r_n.
\end{equation}
Поскольку оценка (\ref{8.09.12-16:01:00}) устанавливается подобно оценке
(\ref{8.09.12-16:00:50}), то мы ограничимся лишь доказательством оценки
(\ref{8.09.12-16:00:50}).

Предположим, что $\alpha (n)\leqslant \frac{1}{p}$. Представим величину $\rho _{n,p+}(f;0)$
в виде суммы трех интегралов
\begin{equation}\label{15.09.12-20:57:48}
    \rho _{n,p+}(f;0)=\int\limits_{|t|\leqslant \alpha (n)n
    }f^\psi _\beta \Big( \frac{t}{n}\Big)\widehat{\tau}_{n,p+}(t)\,dt+$$
    $$+\int\limits_{\alpha (n)n\leqslant |t|\leqslant \frac{n}{p}
    }f^\psi _\beta \Big( \frac{t}{n}\Big)\widehat{\tau}_{n,p+}(t)\,dt+$$
    $$+\int\limits_{|t|\geqslant \frac{n}{p}
    }f^\psi _\beta \Big( \frac{t}{n}\Big)\widehat{\tau}_{n,p+}(t)\,dt.
\end{equation}
Используя для $|t|\leqslant \alpha (n)n$ представление
\begin{equation}\label{11.09.12-21:46:19}
    \widehat{\tau}_{n,p+}(t)=
    \frac{\psi (n)}{\pi }\int\limits_{\frac{n-p}{n}}^{1}
    \frac{nu-n+p}{p}\cos ut\,du+\frac{1}{\pi }
    \int\limits_{1}^{\infty}\psi (nu)\cos ut\,du,
\end{equation}
    получаемое из (\ref{2.12.12-20:32:42}) и (\ref{8.09.12-15:36:59}), находим
\begin{equation}\label{15.09.12-21:07:21}
    \int\limits_{|t|\leqslant \alpha (n)n
    }f^\psi _\beta \Big( \frac{t}{n}\Big)\widehat{\tau}_{n,p+}(t)\,dt=$$
    $$=\frac{\psi (n)}{\pi }\int\limits_{|t|\leqslant \alpha (n)n}f^\psi _\beta \Big( \frac{t}{n}\Big)\int\limits_{
    \frac{n-p}{n}}^{1}\frac{nu-n+p}{p}\cos ut\,du\,dt+
    $$
    $$+\frac{1}{\pi }\int\limits_{|t|\leqslant \alpha (n)n}f^\psi _\beta \Big( \frac{t}{n}\Big)
    \int\limits_{1}^{\infty}
    \psi (nu)\cos ut\,du\,dt.
\end{equation}
Учитывая включение $f^\psi _\beta \in S_\infty$, соотношение
    \begin{equation}\label{15.09.12-22:21:06}
    \Bigg|\int\limits_{\frac{n-p}{n}}^{1}\frac{
    nu-n+p}{p}\cos ut\,du\Bigg|\leqslant \int\limits_{\frac{n-p}{n}}^{1}\frac{
    nu-n+p}{p}\,du=\frac{p}{2n}
    \end{equation}
и условие $\alpha (n)\leqslant 1/p$, из (\ref{15.09.12-21:07:21}) получаем
\begin{equation}\label{15.09.12-21:26:17}
    \int\limits_{|t|\leqslant \alpha (n)n
    }f^\psi _\beta \Big( \frac{t}{n}\Big)\widehat{\tau}_{n,p+}(t)\,dt=$$
    $$=\frac{1}{\pi }\int\limits_{|t|\leqslant \alpha (n)n}f^\psi _\beta \Big( \frac{t}{n}\Big)
    \int\limits_{1}^{\infty}
    \psi (nu)\cos ut\,du\,dt+O(1)\psi (n).
\end{equation}
Как доказано в \cite{Stepanets_VSP} (см. (11.9) и (11.17)), если $f^\psi _\beta \in
S_\infty$, то
\begin{equation}\label{15.09.12-21:43:21}
    \int\limits_{|t|\leqslant \alpha (n)n}f^\psi _\beta \Big( \frac{t}{n}\Big)
    \int\limits_{1}^{\infty}
    \psi (nu)\cos ut\,du\,dt=O(1)\bigg(\psi (n)+\int\limits_{
    1/\alpha (n)}^{\infty}\frac{\psi (t+n)}{t}\,dt\bigg).
\end{equation}
Таким образом, из (\ref{15.09.12-21:26:17}) находим
\begin{equation}\label{15.09.12-21:46:58}
    \int\limits_{|t|\leqslant \alpha (n)n
    }f^\psi _\beta \Big( \frac{t}{n}\Big)\widehat{\tau}_{n,p+}(t)\,dt=
    O(1)\bigg(\psi (n)+\int\limits_{
    1/\alpha (n)}^{\infty}\frac{\psi (t+n)}{t}\,dt\bigg).
\end{equation}

Перейдем к рассмотрению второго интеграла в правой части (\ref{15.09.12-20:57:48}). Для
выделения его главного члена, представим величину $\widehat{\tau}_{n,p+}(t)$ в виде
равенства
\begin{equation}\label{15.09.12-22:04:24}
    \widehat{\tau}_{n,p+}(t)=
    \frac{\psi (n)}{\pi }\int\limits_{\frac{n-p}{n}}^{1}
    \frac{nu-n+p}{p}\cos ut\,du-
    \frac{\psi (n)}{\pi t}\sin t-\frac{n}{\pi t}\int\limits_{
    1}^{\infty}\psi '(nu)\sin ut\,du,
\end{equation}
получаемого из (\ref{11.09.12-21:46:19}) путем преобразования второго интеграла методом
интегрирования по частям:
\begin{equation}\label{15.09.12-13:10:00}
    \int\limits_{1}^{\infty}\psi (nu)\cos ut\,du=-
    \frac{\psi (n)}{t}\sin t-\frac{n}{t}\int\limits_{
    1}^{\infty}\psi '(nu)\sin ut\,du.
\end{equation}
Поскольку, как несложно видеть, исходя из (\ref{15.09.12-22:21:06}),
    $$\int\limits_{\alpha (n)n\leqslant |t|\leqslant
    \frac{n}{p}}f^\psi _\beta \Big( \frac{t}{n}\Big)
    \int\limits_{\frac{n-p}{n}}^{1}\frac{nu-n+p}{p}\cos ut\,du\,dt=O(1),$$
то в силу (\ref{15.09.12-22:04:24}),
\begin{equation}\label{15.09.12-22:26:45}
    \int\limits_{\alpha (n)n\leqslant |t|\leqslant \frac{n}{p}
    }f^\psi _\beta \Big( \frac{t}{n}\Big)\widehat{\tau}_{n,p+}(t)\,dt=$$
    $$=-\frac{\psi (n)}{\pi }\int\limits_{\alpha (n)n\leqslant |t|\leqslant \frac{n
    }{p}}f^\psi _\beta \Big( \frac{t}{n}\Big)\frac{\sin t}{t}\,dt-$$
    $$-\frac{
    n}{\pi }\int\limits_{\alpha (n)n\leqslant |t|\leqslant \frac{n
    }{p}}f^\psi _\beta \Big( \frac{t}{n}\Big)
    \frac{1}{t}\int\limits_{1}^{\infty}\psi '(nu)\sin ut\,du\,dt+O(1)\psi (n).
\end{equation}

Перейдем к третьему интегралу в (\ref{15.09.12-20:57:48}). Дважды интегрируя по частям в
первом интеграле правой части (\ref{15.09.12-22:04:24}), приходим к равенству
\begin{equation}\label{11.09.12-21:50:09}
    \widehat{\tau}_{n,p+}(t)=\frac{\psi (n)}{\pi t^2}\frac{n}{p}\Big(
    \cos t-\cos\frac{n-p}{n}t\Big)-\frac{n}{\pi t}
    \int\limits_{1}^{\infty}\psi '(nu)\sin ut\,du.
\end{equation}
Убеждаясь в справедливости оценки
\begin{equation}\label{15.09.12-22:50:35}
    \frac{n}{p}\int\limits_{|t|\geqslant \frac{n}{p}}f^\psi _\beta \Big( \frac{t}{n}\Big)
    \frac{1}{t^2}\Big(
    \cos t-\cos\frac{n-p}{n}t\Big)\,dt=O(1),
\end{equation}
из (\ref{11.09.12-21:50:09}) и (\ref{15.09.12-22:50:35}) получаем
\begin{equation}\label{15.09.12-22:53:28}
    \int\limits_{|t|\geqslant \frac{n}{p}
    }f^\psi _\beta \Big( \frac{t}{n}\Big)\widehat{\tau}_{n,p+}(t)\,dt=$$
    $$=-\frac{
    n}{\pi }\int\limits_{|t|\geqslant \frac{n
    }{p}}f^\psi _\beta \Big( \frac{t}{n}\Big)
    \frac{1}{t}\int\limits_{1}^{\infty}\psi '(nu)\sin ut\,du\,dt+O(1)\psi (n).
\end{equation}

Объединяя (\ref{15.09.12-21:46:58}), (\ref{15.09.12-22:26:45}), (\ref{15.09.12-22:53:28}) и
оценку
\begin{equation}\label{15.09.12-23:03:00}
    n\int\limits_{|t|\geqslant \alpha (n)n}f^\psi _\beta \Big( \frac{t}{n}\Big)
    \frac{1}{t}\int\limits_{1}^{\infty}\psi '(nu)\sin ut\,du\,dt=$$
    $$=O(1)
    \int\limits_{\alpha (n)}^{\infty}\frac{
    \psi (n)-\psi (n+1/t)}{t}\,dt,
\end{equation}
доказанную в \cite{Stepanets_VSP} (см. (11.31)), из (\ref{15.09.12-20:57:48}) находим
\begin{equation}\label{15.09.12-23:07:51}
    \rho _{n,p+}(f;0)=-\frac{\psi (n)}{\pi }\int\limits_{\alpha (n)n\leqslant |t|\leqslant \frac{n
    }{p}}f^\psi _\beta \Big( \frac{t}{n}\Big)\frac{\sin t}{t}\,dt+$$
    $$+O(1)\bigg(\psi (n)+\int\limits_{
    1/\alpha (n)}^{\infty}\frac{\psi (t+n)}{t}\,dt+
    \int\limits_{\alpha (n)}^{\infty}\frac{
    \psi (n)-\psi (n+1/t)}{t}\,dt\bigg).
\end{equation}

Рассмотрим случай $\alpha (n)>\frac{1}{p}$. Имеет место равенство
\begin{equation}\label{21.09.12-20:57:48}
    \rho _{n,p+}(f;0)=\int\limits_{|t|\leqslant \frac{n}{p}
    }f^\psi _\beta \Big( \frac{t}{n}\Big)\widehat{\tau}_{n,p+}(t)\,dt+$$
    $$+\int\limits_{\frac{n}{p}\leqslant |t|\leqslant \alpha (n)n
    }f^\psi _\beta \Big( \frac{t}{n}\Big)\widehat{\tau}_{n,p+}(t)\,dt+$$
    $$+\int\limits_{|t|\geqslant \alpha (n)n
    }f^\psi _\beta \Big( \frac{t}{n}\Big)\widehat{\tau}_{n,p+}(t)\,dt.
\end{equation}
Воспользовавшись представлением (\ref{11.09.12-21:46:19}) и оценкой (\ref{15.09.12-22:21:06}), получаем
\begin{equation}\label{21.09.12-13:07:00}
    \int\limits_{|t|\leqslant \frac{n}{p}
    }f^\psi _\beta \Big( \frac{t}{n}\Big)\widehat{\tau}_{n,p+}(t)\,dt=$$
    $$=\frac{1}{\pi }\int\limits_{|t|\leqslant \frac{n}{p}}f^\psi _\beta \Big( \frac{t}{n}\Big)
    \int\limits_{1}^{\infty}
    \psi (nu)\cos ut\,du\,dt+O(1)\psi (n).
\end{equation}
В силу (\ref{15.09.12-13:10:00}) и (\ref{11.09.12-21:50:09}),
\begin{equation}\label{21.09.12-13:12:00}
    \widehat{\tau}_{n,p+}(t)=
    \frac{\psi(n)}{\pi}\frac{\sin t}{t}+\frac{1}{\pi}\int\limits_{1}^{\infty}\psi (nu)\cos ut\,du+
    \frac{\psi (n)}{\pi t^2}\frac{n}{p}\Big(
    \cos t-\cos\frac{n-p}{n}t\Big).
\end{equation}
Используя (\ref{21.09.12-13:12:00}) и учитывая, что
    $$\frac{n}{p}\int\limits_{\frac{n}{p}\leqslant |t|\leqslant \alpha (n)n
    }f^\psi _\beta \Big( \frac{t}{n}\Big)\frac{1}{t^2}\Big(
    \cos t-\cos\frac{n-p}{n}t\Big)\,dt=$$
    $$=O(1)\frac{n}{p}\int\limits_{\frac{n}{p}}^{\alpha(n)n}\frac{dt}{t^2}=O(1),$$
получаем для второго интеграла в (\ref{21.09.12-20:57:48})
\begin{equation}\label{21.09.12-13:23:00}
    \int\limits_{\frac{n}{p}\leqslant |t|\leqslant \alpha (n)n
    }f^\psi _\beta \Big( \frac{t}{n}\Big)\widehat{\tau}_{n,p+}(t)\,dt=
    \frac{\psi (n)}{\pi }\int\limits_{\frac{n
    }{p}\leqslant |t|\leqslant \alpha (n)n}f^\psi _\beta \Big( \frac{t}{n}\Big)\frac{\sin t}{t}\,dt+$$
    $$+\frac{1}{\pi }\int\limits_{\frac{n
    }{p}\leqslant |t|\leqslant \alpha (n)n}f^\psi _\beta \Big( \frac{t}{n}\Big)
    \int\limits_{1}^{\infty}
    \psi (nu)\cos ut\,du\,dt+O(1)\psi(n).
\end{equation}
Для нахождения оценки третьего интеграла в правой части (\ref{21.09.12-20:57:48}), воспользуемся представлением (\ref{11.09.12-21:50:09}). Имеем
\begin{equation}\label{25.09.12-12:56:77}
    \int\limits_{|t|\geqslant \alpha (n)n
    }f^\psi _\beta \Big( \frac{t}{n}\Big)\widehat{\tau}_{n,p+}(t)\,dt=$$
    $$=\frac{\psi(n)}{\pi}\frac{n}{p}\int\limits_{|t|\geqslant \alpha (n)n
    }f^\psi _\beta \Big( \frac{t}{n}\Big)\frac{1}{t^2}\Big(\cos t-\cos \frac{n-p}{n}t\Big)\,dt-$$
    $$-\frac{n}{\pi}\int\limits_{|t|\geqslant \alpha (n)n
    }f^\psi _\beta \Big( \frac{t}{n}\Big)\frac{1}{t}\int\limits_{1}^{\infty}\psi'(nu)\sin ut\,du\,dt.
\end{equation}
Поскольку, как несложно видеть
    $$\frac{n}{p}\int\limits_{|t|\geqslant \alpha (n)n
    }f^\psi _\beta \Big( \frac{t}{n}\Big)\frac{1}{t^2}\Big(\cos t-\cos \frac{n-p}{n}t\Big)\,dt=
    O(1)\frac{n}{p}\int\limits_{\alpha(n)n}^{\infty}\frac{dt}{t^2}=O(1),$$
то, в силу (\ref{15.09.12-23:03:00}), из (\ref{25.09.12-12:56:77}) получаем
\begin{equation}\label{25.09.12-12:56:00}
    \int\limits_{|t|\geqslant \alpha (n)n
    }f^\psi _\beta \Big( \frac{t}{n}\Big)\widehat{\tau}_{n,p+}(t)\,dt=O(1)\Bigg(\psi(n)+
    \int\limits_{\alpha (n)}^{\infty}\frac{
    \psi (n)-\psi (n+1/t)}{t}\,dt\Bigg).
\end{equation}
Объединяя (\ref{21.09.12-20:57:48}), (\ref{21.09.12-13:07:00}), (\ref{21.09.12-13:23:00}),
(\ref{25.09.12-12:56:00}) и учитывая (\ref{15.09.12-21:43:21}), приходим к оценке
\begin{equation}\label{25.09.12-13:04:51}
    \rho _{n,p+}(f;0)=\frac{\psi (n)}{\pi }\int\limits_{\frac{n
    }{p}\leqslant |t|\leqslant \alpha (n)n}f^\psi _\beta \Big( \frac{t}{n}\Big)\frac{\sin t}{t}\,dt+$$
    $$+O(1)\bigg(\psi (n)+\int\limits_{
    1/\alpha (n)}^{\infty}\frac{\psi (t+n)}{t}\,dt+
    \int\limits_{\alpha (n)}^{\infty}\frac{
    \psi (n)-\psi (n+1/t)}{t}\,dt\bigg).
\end{equation}
Из (\ref{15.09.12-23:07:51}) и (\ref{25.09.12-13:04:51}) следует справедливость
(\ref{8.09.12-16:00:50}) при $\alpha(n) >\frac{1}{p}$.

Для доказательства леммы в случае, когда $\alpha (n)>1$, достаточно применить к первому
интегралу в (\ref{25.09.12-13:04:51}) установленную в работе \mbox{\cite[c.
119]{Stepanets_1987_Izv}} (лемма 5) оценку
\begin{equation}\label{12.11.12-22:55:54}
    \Bigg|\int\limits_{|t|\geqslant \alpha^* (n)}\varphi \big(
    \frac{t}{n}\big)\frac{\sin\Big(t+\frac{\gamma  \pi }{2}\Big)}{t}\,dt\Bigg|=O(1),\ \
    \varphi \in S_\infty,\ \ n\to\infty,\ \ \gamma \in\mathbb{R},
\end{equation}
где $\{\alpha ^*(n)\}=\{\alpha (n): \alpha (n)\geqslant n\pi \}.$ Действительно, используя
(\ref{12.11.12-22:55:54}), получаем
\begin{equation}\label{12.11.12-23:15:03}
    \int\limits_{\frac{n
    }{p}\leqslant |t|\leqslant \alpha (n)n}f^\psi _\beta \Big( \frac{t}{n}\Big)\frac{\sin t}{t}\,dt
    =\int\limits_{\frac{n
    }{p}\leqslant |t|\leqslant n}
    f^\psi _\beta \Big( \frac{t}{n}\Big)\frac{\sin t}{t}\,dt+O(1).
\end{equation}
Из (\ref{25.09.12-13:04:51}) и (\ref{12.11.12-23:15:03}) следует оценка
(\ref{8.09.12-16:00:50}) в случае, когда $\alpha(n)>1$. Лемма \ref{8.09.12-15:54:03}
доказана.
\end{proof}

Положив в условии леммы \ref{8.09.12-15:54:03}
    $$\alpha(n)=\alpha(\psi;n)=\frac{1}{\eta(\psi;n)-n},$$
где $\eta(\psi;n)$ определяется формулой (\ref{27.11.12-23:55:28}), а $\psi \in F$
($F\subset\mathfrak{M}'$) и учитывая доказанные в работе \cite[c. 155, 156]{Stepanets_VSP}
оценки
    $$\int\limits_{\eta(\psi;n)}^\infty\frac{\psi(t)}{t-n}\,dt=O(1)\psi(n),\ \ \forall \psi\in F,$$
    $$\int\limits_n^{\eta(\psi;n)}\frac{\psi(n)-\psi(t)}{t-n}\, dt=O(1)\psi(n),\ \ \forall \psi\in F$$
и
    \begin{equation}\label{25.11.12-01:43:57}
    \mu (n)=\mu (\psi ;n)=\frac{n}{\eta(\psi ;n)-n}\geqslant K>0,\ \ \forall \psi\in F,
    \end{equation}
получаем следующее утверждение.

\begin{cor} Пусть $\psi\in F$, $\beta\in \mathbb{R}$, $n,p\in \mathbb{N},$ $p\leqslant n.$ Тогда для произвольной функции $f\in C^\psi_{\beta,
\infty}$ при $n\to\infty$ справедливо равенство
\begin{equation}\label{10.10.12-13:21:00}
    \rho _{n,p}(f;0)=(-1)^s\frac{\psi (n)}{\pi }\int\limits_{\mathcal{I}_{n}}
    f^\psi _\beta \Big( \frac{t}{n}\Big)\frac{\sin \Big(t+\frac{\beta\pi}{2}\Big)}{t}\,dt+O(1)\psi(n),
\end{equation}
где
    $$\mathcal{I}_{n}=\Big(-\frac{n}{p},-\mu(n)\Big)\bigcup
    \Big(\mu(n), \frac{n}{p}\Big),\ \ \text{если $\mu(n)\leqslant \frac{n}{p}$},$$
    $$\mathcal{I}_{n}=\Big(-\mu(n), -\frac{n}{p}\Big)\bigcup
    \Big(\frac{n}{p}, \mu(n)\Big),\ \ \text{если $\frac{n}{p}<\mu (n)\leqslant n$},$$
    $$\mathcal{I}_{n}=\Big(-n, -\frac{n}{p}\Big)\bigcup
    \Big(\frac{n}{p}, n\Big),\ \ \text{если $\mu (n)> n$},$$
    $$\mu(n)=\mu (\psi ;n)=\frac{n}{\eta(\psi;n)-n},$$
\begin{equation}\label{7.11.12-21:36:30}
    s=s(\psi, n, p)=
  \begin{cases}
    1, & \mu (n)\leqslant \frac{n}{p}, \\
    0, & \mu (n)> \frac{n}{p},
  \end{cases}
\end{equation} а $O(1)$
--- величина, равномерно ограниченная по $\beta,$ $n$ и $p.$
\end{cor}

Далее, положим
    $$t_k=(2k+1-\beta )\frac{\pi }{2},\ \ \ k\in \mathbb{Z},$$
    $$x_k=t_k+\frac{\pi }{2}, \ \ \ k\in \mathbb{Z}$$
и
\begin{equation}\label{13.11.12-22:01:24}
l_n(t)=
  \begin{cases}
    \frac{1}{x_k}, & t\in[t_{k}, t_{k+1}],\ \ k=k_0, k_{0}+1,\ldots k_1-1,k_2,k_2+1,\ldots k_3-1,\\
    0, & t\in(-\infty, t_{k_0})\cup(t_{k_3}, \infty),
  \end{cases}
\end{equation}
где $k_0,$ $k_1,$ $k_2$ и $k_3$ выбраны таким образом, чтобы
    $$t_{k_0-1}<-\frac{n}{p}\leqslant t_{k_0},$$
    $$t_{k_1}<-\mu (n)\leqslant t_{k_1+1},$$
    $$t_{k_2-1}<\mu (n)\leqslant t_{k_2},$$
    $$t_{k_3}<\frac{n}{p}\leqslant t_{k_3+1},$$
если $\mu(n)\leqslant \frac{n}{p},$
    $$t_{k_0-1}<-\mu (n)\leqslant t_{k_0},$$
    $$t_{k_1}<-\frac{n}{p}\leqslant t_{k_1+1},$$
    $$t_{k_2-1}<\frac{n}{p}\leqslant t_{k_2},$$
    $$t_{k_3}<\mu (n)\leqslant t_{k_3+1},$$
если $\frac{n}{p}< \mu(n)\leqslant n,$ и
$$t_{k_0-1}<-n\leqslant t_{k_0},$$
    $$t_{k_1}<-\frac{n}{p}\leqslant t_{k_1+1},$$
    $$t_{k_2-1}<\frac{n}{p}\leqslant t_{k_2},$$
    $$t_{k_3}<n\leqslant t_{k_3+1},$$
если $\mu(n)> n.$

С учетом принятых обозначений, интеграл в правой части равенства (\ref{10.10.12-13:21:00})
можно представить в виде
\begin{equation}\label{20.10.12-19:54:23}
    \int\limits_{\mathcal{I}_{n}}
    f^\psi _\beta \Big( \frac{t}{n}\Big)\frac{\sin \Big(t+\frac{\beta\pi}{2}\Big)}{t}\,dt=
    \int\limits_{
    (t_{k_0}, t_{k_1})\cup(t_{k_2}, t_{k_3})}f^\psi _\beta \Big( \frac{t}{n}\Big)l_n(t)
    \sin \Big(t+\frac{\beta\pi}{2}\Big)\,dt+$$
    $$+R_{n,p}^{(1)}+R_{n,p}^{(2)},
\end{equation}
где
    \begin{equation}\label{27.10.12-18:31:42}
    R_{n,p}^{(1)}=\int\limits_{
    (t_{k_0}, t_{k_1})\cup(t_{k_2}, t_{k_3})}f^\psi _\beta \Big( \frac{t}{n}\Big)\Big(\frac{
    1}{t}-l_n(t)\Big)
    \sin \Big(t+\frac{\beta\pi}{2}\Big)\,dt,
    \end{equation}
    \begin{equation}\label{27.10.12-18:31:52}
    R_{n,p}^{(2)}=
    \bigg(\int\limits_{-a}^{t_{k_0}}+\int\limits_{t_{k_1}}^{-b}+
    \int\limits_{b}^{t_{k_2}}
    +\int\limits_{t_{k_3}}^{a}\bigg)
    f^\psi _\beta \Big( \frac{t}{n}\Big)\frac{\sin \Big(t+\frac{\beta\pi}{2}\Big)}{t}\,dt,
    \end{equation}
    \begin{equation}\label{8.11.12-23:44:31}
    a=
\begin{cases}
    \frac{n}{p}, & \mu (n)\leqslant \frac{n}{p}, \\
    \mu (n), & \frac{n}{p}<\mu (n)\leqslant n,\\
    n, & \mu (n)>n,
  \end{cases}
\end{equation}
\begin{equation}\label{2.12.12-20:37:31}
    b=
  \begin{cases}
    \mu (n), & \mu (n)\leqslant \frac{n}{p}, \\
    \frac{n}{p}, & \mu (n)>\frac{n}{p}.
  \end{cases}
\end{equation}

Покажем, что имеют место оценки
\begin{equation}\label{20.10.12-20:30:25}
    R_{n,p}^{(1)}=O(1)
\end{equation}
и
\begin{equation}\label{20.10.12-20:30:58}
    R_{n,p}^{(2)}=O(1).
\end{equation}
Выполняя элементарные преобразования в (\ref{27.10.12-18:31:42}) и учитывая, что $|f^\psi
_\beta (\cdot)|\leqslant 1$, имеем
 \begin{equation}\label{20.10.12-22:45:02}
 R_{n,p}^{(1)}=O(1)\bigg(\sum\limits_{k=k_0}^{k_1-1}\int\limits_{
 t_k}^{t_{k+1}}\Big|\frac{1}{t}-\frac{1}{x_k}\Big|\,dt+
 \sum\limits_{k=k_2}^{k_3-1}\int\limits_{
 t_k}^{t_{k+1}}\Big|\frac{1}{t}-\frac{1}{x_k}\Big|\,dt\bigg).
 \end{equation}
Поскольку
    $$\int\limits_{
 t_k}^{x_{k}}\Big|\frac{1}{t}-\frac{1}{x_k}\Big|\,dt>\int\limits_{
 x_k}^{t_{k+1}}\Big|\frac{1}{t}-\frac{1}{x_k}\Big|\,dt,\ \ k\in\mathbb{Z},$$
то
    \begin{equation}\label{20.10.12-22:57:40}
    \int\limits_{
 t_k}^{t_{k+1}}\Big|\frac{1}{t}-\frac{1}{x_k}\Big|\,dt<2\int\limits_{
 t_k}^{x_{k}}\Big|\frac{1}{t}-\frac{1}{x_k}\Big|\,dt=2\int\limits_{
 t_k}^{x_{k}}\frac{x_k-t}{|tx_k|}\,dt\leqslant \pi \int\limits_{
 t_k}^{x_{k}}\frac{dt}{t^2}<\pi \int\limits_{t_k}^{t_{k+1}}\frac{dt}{t^2}.
    \end{equation}
Объединяя (\ref{20.10.12-22:45:02}) и (\ref{20.10.12-22:57:40}), получаем
\begin{equation}\label{20.10.12-23:00:23}
    R_{n,p}^{(1)}=O(1)\bigg(\int\limits_{t_{k_0}}^{t_{k_1}}\frac{dt}{t^2}+
    \int\limits_{t_{k_2}}^{t_{k_3}}\frac{dt}{t^2}\bigg)=
    O(1)\bigg(\frac{1}{|t_{k_1}|}+\frac{1}{t_{k_2}}\bigg)=
    O(1)\max\Big\{\frac{1}{\mu (n)},1 \Big\}.
\end{equation}
Поскольку в силу (\ref{25.11.12-01:43:57}) $\mu (n)$  ограничена снизу, то из
(\ref{20.10.12-23:00:23}) следует оценка (\ref{20.10.12-20:30:25}).

Учитывая соотношения
    $$t_{k_0}+a<t_{k_0}-t_{k_0-1}=\pi, $$
    $$-b-t_{k_1}\leqslant t_{k_1+1}-t_{k_1}=\pi ,$$
    $$t_{k_2}-b<t_{k_2}-t_{k_{2}-1}=\pi $$
и
    $$a-t_{k_3}\leqslant t_{k_3+1}-t_{k_3}=\pi ,$$
после несложных преобразований из (\ref{27.10.12-18:31:52}) находим
    $$R_{n,p}^{(2)}=O(1)\bigg(
    \frac{t_{k_0}+a}{|t_{k_0}|}+\frac{-b-t_{k_1}}{b}+\frac{t_{k_2}-b}{b}+
    \frac{a-t_{k_3}}{t_{k_3}}
    \bigg)=$$
    $$=O(1)\bigg(
    \frac{1}{|t_{k_0}|}+\frac{1}{b}+
    \frac{1}{t_{k_3}}
    \bigg)=O(1)\max\Big\{\frac{1}{\mu (n)},1 \Big\}=
    O(1).$$

Объединяя (\ref{10.10.12-13:21:00}), (\ref{20.10.12-19:54:23}), (\ref{20.10.12-20:30:25}) и
(\ref{20.10.12-20:30:58}), приходим к формуле
\begin{equation}\label{27.10.12-19:16:56}
    \rho _{n,p}(f;0)=(-1)^s\frac{\psi (n)}{\pi }\int\limits_{
    (t_{k_0}, t_{k_1})\cup(t_{k_2}, t_{k_3})}f^\psi _\beta \Big( \frac{t}{n}\Big)l_n(t)
    \sin \Big(t+\frac{\beta\pi}{2}\Big)\,dt+$$
    $$+O(1)\psi (n).
\end{equation}

Приступим к нахождению оценки величины $\mathcal{E}(C^\psi _{\beta ,\infty};U_{n,p}^\psi)$.
Подставляя (\ref{27.10.12-19:16:56}) в (\ref{1.12.12-16:19:29}), имеем
\begin{equation}\label{27.10.12-22:49:31}
    \mathcal{E}(C^\psi _{\beta ,\infty};U_{n,p}^\psi)\leqslant
    \frac{\psi (n)}{\pi }
    \int\limits_{
    (t_{k_0}, t_{k_1})\cup(t_{k_2}, t_{k_3})}\Big|l_n(t)\sin\Big(t+\frac{\beta \pi }{2}
    \Big)\Big|\,dt+O(1)\psi (n).
\end{equation}
Покажем, что
\begin{equation}\label{27.10.12-22:52:10}
    \int\limits_{
    (t_{k_0}, t_{k_1})\cup(t_{k_2}, t_{k_3})}\Big|l_n(t)\sin\Big(t+\frac{\beta \pi }{2}
    \Big)\Big|\,dt=\frac{4}{\pi }A_{n,p}^\psi +O(1),
\end{equation}
где
\begin{equation}\label{13.11.12-00:08:21}
    A_{n,p}^\psi =
  \begin{cases}
    \ln\frac{n}{p\mu (n)}, & \text{если}\,\,\,\, \mu (n)\leqslant  \frac{n}{p}, \\
    &\\
    \ln \frac{p\mu (n)}{n}, & \text{если}\,\,\,\, \frac{n}{p}<\mu (n)\leqslant n,\\
    &\\
    \ln p, & \text{если}\,\,\,\, \mu (n)>n.
  \end{cases}
\end{equation}

Учитывая центрально-симметричность относительно точки $x_k$ на $[t_k, t_{k+1}],$ $k\in
\mathbb{Z}$ функции \mbox{$\sin(t+\beta \pi /2)$}, имеем
    $$\int\limits_{t_k}^{t_{k+1}}\Big|\sin\Big(t+\frac{\beta \pi }{2}
    \Big)\Big|\,dt=2\int\limits_{t_k}^{x_{k}}\Big|\sin\Big(t+\frac{\beta \pi }{2}
    \Big)\Big|\,dt=$$
    $$=2\int\limits_{(k+\frac{1}{2})\pi }^{(k+1)\pi }|\sin t|\,dt=2\int\limits_{
    0}^{\frac{\pi }{2}}\cos t\,dt=2.$$
Отсюда получаем, что
\begin{equation}\label{27.10.12-23:01:10}
    \int\limits_{
    (t_{k_0}, t_{k_1})\cup(t_{k_2}, t_{k_3})}\Big|l_n(t)\sin\Big(t+\frac{\beta \pi }{2}
    \Big)\Big|\,dt=\sum\limits_{
    k=k_0}^{k_1-1}\frac{1}{|x_k|}\int\limits_{t_k}^{t_{k+1}}\Big|\sin\Big(t+\frac{\beta \pi }{2}
    \Big)\Big|\,dt+$$
    $$+\sum\limits_{
    k=k_2}^{k_3-1}\frac{1}{x_k}\int\limits_{t_k}^{t_{k+1}}\Big|\sin\Big(t+\frac{\beta \pi }{2}
    \Big)\Big|\,dt=2\bigg(\sum\limits_{
    k=k_0}^{k_1-1}\frac{1}{|x_k|}+\sum\limits_{
    k=k_2}^{k_3-1}\frac{1}{x_k}\bigg).
\end{equation}
Поскольку
    $$t_{k+1}-t_k=\pi ,\ \ k\in\mathbb{Z},$$
то имеет место равенство
    \begin{equation}\label{9.11.12-22:44:23}
    \sum\limits_{k=k_0}^{k_1-1}\frac{1}{|x_k|}=
    \frac{1}{\pi }\sum\limits_{k=k_0}^{k_1-1}\int\limits_{t_k}^{t_{k+1}}\frac{1}{|x_k|}\,dt=
    -\frac{1}{\pi }\int\limits_{t_{k_0}}^{t_{k_1}}\frac{dt}{t}+
    \frac{1}{\pi }\sum\limits_{k=k_0}^{k_1-1}\int\limits_{t_k}^{t_{k+1}}\bigg(
    -\frac{1}{x_k}+\frac{1}{t}\bigg)\,dt.
    \end{equation}
Используя оценку (\ref{20.10.12-22:57:40}), получаем
\begin{equation}\label{9.11.12-22:36:37}
    \bigg|\sum\limits_{k=k_0}^{k_1-1}\int\limits_{t_k}^{t_{k+1}}\bigg(
    -\frac{1}{x_k}+\frac{1}{t}\bigg)\,dt\bigg|\leqslant
    \sum\limits_{k=k_0}^{k_1-1}\int\limits_{t_k}^{t_{k+1}}\bigg|
    \frac{1}{t}-\frac{1}{x_k}\bigg|\,dt<$$
    $$<\pi \sum\limits_{k=k_0}^{k_1-1}\int\limits_{
    t_k}^{t_{k+1}}\frac{dt}{t^2}=\pi \int\limits_{
    t_{k_0}}^{t_{k_1}}\frac{dt}{t^2}< \frac{\pi }{|t_{k_1}|}\leqslant K \max\Big\{
    \frac{1}{\mu (n)}, 1\Big\}=O(1).
\end{equation}
Сопоставляя (\ref{9.11.12-22:44:23}) и (\ref{9.11.12-22:36:37}), находим
\begin{equation}\label{9.11.12-22:46:45}
    \sum\limits_{k=k_0}^{k_1-1}\frac{1}{|x_k|}=
    -\frac{1}{\pi }\int\limits_{t_{k_0}}^{t_{k_1}}\frac{dt}{t}+O(1).
\end{equation}
Используя обозначения (\ref{8.11.12-23:44:31}) и (\ref{2.12.12-20:37:31}), запишем
(\ref{9.11.12-22:46:45}) в виде
\begin{equation}\label{9.11.12-22:50:19}
    \sum\limits_{k=k_0}^{k_1-1}\frac{1}{|x_k|}=-\frac{1}{\pi }\int\limits_{
    -a}^{-b}\frac{dt}{t}+\frac{1}{\pi }\int\limits_{
    -a}^{-b}\frac{dt}{t}-\frac{1}{\pi }\int\limits_{t_{k_0}}^{t_{k_1}}\frac{dt}{t}+O(1)=$$
    $$=-\frac{1}{\pi }\int\limits_{
    -a}^{-b}\frac{dt}{t}+\int\limits_{t_{k_1}}^{-b}\frac{dt}{t}+
    \int\limits_{-a}^{t_{k_0}}\frac{dt}{t}+O(1).
\end{equation}
Поскольку
    $$-b-t_{k_1}\leqslant \pi $$
и
    $$t_{k_0}+a<\pi ,$$
то
    $$\int\limits_{t_{k_1}}^{-b}\frac{dt}{t}+
    \int\limits_{-a}^{t_{k_0}}\frac{dt}{t}=O(1)\Big(\frac{1}{b}+\frac{1}{|t_{k_0}|}\Big)
    =O(1)\frac{1}{|t_{k_0}|}=$$
\begin{equation}\label{9.11.12-23:05:20}
    =O(1)\max\Big\{
    \frac{1}{\mu (n)}, 1\Big\}=O(1).
\end{equation}
Объединяя (\ref{9.11.12-22:50:19}) и (\ref{9.11.12-23:05:20}), получаем
$$
    \sum\limits_{k=k_0}^{k_1-1}\frac{1}{|x_k|}=-\frac{1}{\pi }\int\limits_{
    -a}^{-b}\frac{dt}{t}+O(1)=\frac{1}{\pi }\ln \frac{a}{b}+O(1).
$$
Таким образом
    \begin{equation}\label{9.11.12-23:07:00}
    \sum\limits_{k=k_0}^{k_1-1}\frac{1}{|x_k|}=\frac{1}{\pi }A_{n,p}^\psi +O(1),
    \end{equation}
где $A_{n,p}$ определяется формулой (\ref{13.11.12-00:08:21}).

Аналогичным образом доказывается оценка
\begin{equation}\label{9.11.12-23:10:58}
    \sum\limits_{k=k_2}^{k_3-1}\frac{1}{x_k}=\frac{1}{\pi }A_{n,p}^\psi +O(1).
\end{equation}
Из (\ref{27.10.12-23:01:10}), (\ref{9.11.12-23:07:00}) и (\ref{9.11.12-23:10:58}) следует
соотношение (\ref{27.10.12-22:52:10}).

Сопоставляя (\ref{27.10.12-22:49:31}) и (\ref{27.10.12-22:52:10}), находим
\begin{equation}\label{9.11.12-23:17:02}
    \mathcal{E}(C^\psi _{\beta ,\infty};U_{n,p}^\psi)\leqslant
    \psi (n)
    \bigg(\frac{4}{\pi^2 }A_{n,p}^\psi +O(1)\bigg).
\end{equation}
Для окончательного доказательства теоремы остается показать, что в (\ref{9.11.12-23:17:02})
можно поставить знак равенства.

Обозначим через $\varphi _0(t)$ такую $2\pi $-периодическую функцию, что
    $$\int\limits_{-\pi }^{\pi }\varphi _0(t)\,dt=0,\ \ \ |\varphi _0(t)|\leqslant 1$$
и
    $$\varphi _0(t)=\text{sign}\bigg(l_n(nt)\sin\Big(nt+\frac{\beta \pi }{2}\Big)\bigg), \
    \ \ t\in[-1,1],$$
где $l_n(t)$ определяется согласно (\ref{13.11.12-22:01:24}). Очевидно, что такая функция
существует. Согласно п. 7.2 работы \cite[c. 136, 137]{Stepanets_1995}, в классе
$C^{\psi}_{\beta,\infty},$ $\psi \in F$ найдется функция $f_0(t),$ для которой $\varphi
_0(t)$ будет ее $(\psi ,\beta )$-производной. Для $f_0(t)$ из (\ref{27.10.12-19:16:56}) и
(\ref{27.10.12-22:52:10}) получаем равенство
\begin{equation}\label{13.11.12-21:47:09}
    |\rho _{n,p}(f_0;0)|=\frac{\psi (n)}{\pi }\bigg|\int\limits_{
    (t_{k_0}, t_{k_1})\cup(t_{k_2}, t_{k_3})}\varphi_0 \Big( \frac{t}{n}\Big)l_n(t)
    \sin \Big(t+\frac{\beta\pi}{2}\Big)\,dt\bigg|+$$
    $$+O(1)\psi (n)=\frac{\psi (n)}{\pi }\int\limits_{
    (t_{k_0}, t_{k_1})\cup(t_{k_2}, t_{k_3})}\bigg|l_n(t)
    \sin \Big(t+\frac{\beta\pi}{2}\Big)\bigg|\,dt+O(1)\psi (n)=$$
    $$=\psi (n)
    \bigg(\frac{4}{\pi^2 }A_{n,p}^\psi +O(1)\bigg).
\end{equation}
Поскольку $\mathcal{E}(C^\psi _{\beta ,\infty};U_{n,p}^\psi)\geqslant |\rho _{n,p}(f_0,0)|$,
то из (\ref{9.11.12-23:17:02}) и (\ref{13.11.12-21:47:09}) вытекает справедливость теоремы
\ref{17.11.12-11:21:14}.

\end{document}